\newcommand{\ba}{\begin{array}}
\newcommand{\ea}{\end{array}}
\newcommand{\bc}{\begin{center}}
\newcommand{\ec}{\end{center}}

\newcommand{\beqn}[1]{\begin{equation}\label{#1}}
\newcommand{\eeqn}{\end{equation}}
\newcommand{\be}{\begin{equation}}
\newcommand{\ee}{\end{equation}}

\newcommand{\beqnn}{\begin{eqnarray}}
\newcommand{\eeqnn}{\end{eqnarray}}

\newtheorem{theorem}{Theorem}
\newtheorem{remark}{Remark}

\newtheorem{lemma}{Lemma}
\newtheorem{definition}{Definition}

\newtheorem{assumption}{Assumption}
\newcommand{\T}{{\rm T}}
\documentclass{IEEEtran}
\usepackage{cite}
\usepackage{amsmath,amssymb,amsfonts}
\usepackage{algorithmic}
\usepackage{graphicx}
\usepackage{textcomp}
\usepackage{cases}
\usepackage{eurosym}
\usepackage{booktabs}
\usepackage{multirow}
\usepackage{algorithm}
\usepackage{algorithmic}
\usepackage{xcolor}

\def\BibTeX{{\rm B\kern-.05em{\sc i\kern-.025em b}\kern-.08em
    T\kern-.1667em\lower.7ex\hbox{E}\kern-.125emX}}
\begin{document}
\title{Self-Triggered Adaptive Model Predictive Control of Constrained Nonlinear Systems: A Min-Max Approach}
\author{Kunwu Zhang, \IEEEmembership{Student Member,~IEEE,} Changxin Liu, \IEEEmembership{Student Member,~IEEE,} and Yang Shi, \IEEEmembership{Fellow,~IEEE}
\thanks{The authors are with the Department of Mechanical Engineering, University of Victoria, Victoria, BC, V8W 2Y2, Canada (e-mail: kunwu@uvic.ca; chxliu@uvic.ca; yshi@uvic.ca).}
}
\maketitle

\begin{abstract}
In this paper, a self-triggered adaptive  model predictive control (MPC) algorithm is proposed for constrained discrete-time nonlinear systems subject to parametric uncertainties and additive disturbances. To bound the parametric uncertainties with reduced overestimation, a zonotope-based set-membership parameter estimator is developed, which is also compatible with the aperiodic sampling resulted from the self-triggering mechanism. The estimation of uncertainties is employed to reformulate the optimization problem in a min-max MPC scheme to reduce the conservatism. By designing a time-varying penalty in the cost function, the estimation of uncertainties is implicitly considered in the self-triggering scheduler, therefore making the triggering interval further optimized. The resulting self-triggered adaptive MPC algorithm guarantees the recursive feasibility, while providing less conservative performance compared with the self-triggered robust MPC method. Furthermore, we theoretically show that the closed-loop system is input-to-state practical stable (ISpS) at triggering time instants. A numerical example and comparison study are performed to demonstrate the efficacy of the proposed method.
\end{abstract}


\section{Introduction}\label{STAMPC:introduction}
The problem of addressing the computation and communication constraints explicitly in networked dynamic systems has attracted increasing attention in recent years \cite{heemels_introduction_2012CDC}. Compared with periodic implementations, the event-based aperiodic control is a more promising solution to achieve the trade-off between the closed-loop performance and the overall communication load, since the control input is not calculated and transmitted until a certain well-defined event related to closed-loop behaviors occurs. Such an event is generally triggered at time instants when the system output or state leaves a certain set \cite{heemels_introduction_2012CDC,tabuada_event_2007TAC,wang_event_2010TAC}. Hence, even-triggered control requires continuously monitoring system states to determine the computation and communication of control inputs, which may be infeasible for some networked systems with limited communication resources. To further reduce the communication load, the self-triggered approaches have been proposed, where the next sampling time instant is determined by the triggering scheduler at the current time instant so that the system states are only measured at triggering time instants \cite{velasco_self_RTSS2003}. %But this advantage may render the system more susceptible to disturbances and uncertainties. 
A comprehensive introduction to event- and self-triggered control can be referred to \cite{heemels_introduction_2012CDC}.

In past decades, model predictive control (MPC) has achieved the phenomenal success in process industries due to its capability of efficiently handling hard constraints on inputs and states for complicated systems \cite{qin_survey_2003CEP}. In MPC, the control input is obtained by solving a finite-horizon optimal control problem at each time instant \cite{acipl_li_robust_2016,mayne_mpcsurvey_2014auto}. This strategy ensures the optimal performance with respect to a certain performance index, however, it inherently introduces the increased computational complexity and hence may restrict its application to many practical control problems. Unlike conventional periodic MPC, in event-triggered MPC the new control inputs are only computed and transmitted if a certain triggering threshold is 
reached, consequently saving communication and computational power,  e.g. \cite{acipl_li_event_2014automatica,acipl_li_periodic_2015SCL,acipl_liu_codeisign_2019TAC}.

Compared with the event-triggered method, using the self-triggered approach can not only reduce the average frequency of computing the control input, but also reduce the overall communication load since the system state or output is only measured and transmitted at triggering time instants\cite{heemels_introduction_2012CDC}. Some results addressing self-triggered MPC have been reported in the literature, e.g., \cite{hashimoto_self_2016TAC,acipl_Chen_2018IS,sun_robust_2019TAC,acipl_li_triggering_2018TAC,acipl_liu_robust_2019Auto,brunner_robust_2014ECC,brunner_robust_2016Automatica}. A constrained nonlinear system is considered in \cite{acipl_li_triggering_2018TAC}, where the authors proposed a co-design strategy such that the maximum triggering interval and the
optimal control inputs can be simultaneously obtained by solving an optimization problem. For disturbed linear systems, a tube-based self-triggered MPC algorithm is presented in \cite{brunner_robust_2014ECC}, where the static state tube is constructed to guarantee the robust constraint satisfaction. To enlarge the region of attraction, the integration of self-triggered MPC with homothetic tubes is proposed in \cite{brunner_robust_2016Automatica}, where the state tubes are optimized online to reduce the conservatism.
%additional scaling factors are introduced as decision variables for MPC optimization problems so that the state tubes are optimized online. 
The probabilistic constraints and stochastic disturbances are considered in  \cite{dai_stochastic_2018Automatica,acipl_Chen_2018IS}.
%In , the solutions handling stochastic disturbances and constraints for linear systems are developed by synthesizing the self-triggering mechanism and stochastic MPC. 
For nonlinear systems, in \cite{hashimoto_self_2016TAC} a robust self-triggered MPC scheme is developed for unconstrained nonlinear affine systems, where the triggering interval is maximized by comparing the current optimal cost and the predicted optimal cost. To relieve the computational burden, the authors in \cite{sun_robust_2019TAC} have proposed an adaptive mechanism for the prediction horizon in the dual mode MPC framework. By combining the self-triggering mechanism with the min-max MPC strategy, a recent work in \cite{acipl_liu_robust_2019Auto} provided a novel robust self-triggered MPC algorithm for general nonlinear systems considering both parametric uncertainties and additive disturbances. 

In this work, we investigate self-triggered adaptive MPC for discrete-time nonlinear systems subject to both parametric uncertainties and additive disturbances. Our primary motivation is as follows. In most of aforementioned works on self-triggered robust MPC approaches, the uncertainty is handled by considering its worst-case realization or tightening state constraints, which essentially and heavily relies on the initial guess for bounds on uncertainties. Inherently, those methods are conservative for handling fixed or slowly changing uncertainties. It has been found in studies that adaptive MPC is a promising solution to mitigate the conservatism of robust MPC \cite{adetola_adaptiveMPC_2009SCL,lorenzen_robust_2019Automatica}, where its main insight is to accommodate the online estimation of uncertainty within a robust MPC framework. To the best of our knowledge, self-triggered adaptive MPC for nonlinear systems has not been studied since synthesizing uncertainty estimation with self-triggered robust MPC introduces some new theoretical and practical problems. For example, estimating the uncertainty based on the input and state history may result in recursive updates of the system model, deteriorating the performance or even destroy the closed-loop stability, especially for nonlinear systems. In addition, the self-triggering scheduler makes the system sampled aperiodically, which renders the uncertainty estimation scheme infeasible. Another remarkable difficulty lies in guaranteeing the robust constraint satisfaction for the aperiodically sampled system with online changing models.
%Hence, a suitable uncertainty estimation routine adapted to the aperiodic sampling is desired to obtain an accurate description of uncertainty. 

To solve these problems, we develop a self-triggered adaptive  MPC algorithm with the following features: 1) A zonotope-based set-membership parameter estimator is developed to approximate the feasible solution set (FSS) of unknown parameters with reduced overestimation. By estimating the reachable set of system states, the proposed estimator becomes compatible with the aperiodic sampling. 2) The co-design of MPC optimization and triggering time instants is considered in the proposed self-triggering mechanism. According to the estimated FSS (EFSS), we firstly construct the cost function consisting of the penalized stage costs with the open loop prediction, the stage costs with the closed-loop prediction and the terminal cost. The open loop stage costs are penalized with time-varying weights so that the new estimation of uncertainty is implicitly considered in the proposed self-triggering mechanism. By comparing the optimal cost with different open loop scenarios, the optimal triggering intervals are determined, therefore leading to the reduced average sampling frequency in the closed-loop system. 3) A self-triggered adaptive MPC algorithm is proposed based on the min-max MPC technique.  We effectively facilitate the online parameter adaption in the proposed MPC scheme by reformulating the optimization problem based on the new EFSS. The resulting self-triggered adaptive MPC algorithm guarantees the recursive feasibility, while providing comparable closed-loop performance and reduced average sampling frequency compared with the self-triggered robust MPC method. It is also theoretically shown that the closed-loop system is input-to-state practical stable (ISpS) at triggering time instants. 

The remainder of this paper is organized as follows: In Section \ref{STAMPC:sec:problem}, the problem setup is demonstrated. Section \ref{STAMPC:sec:estimation} describes the design of the set-membership parameter estimator under the self-triggering mechanism. An self-triggered adaptive MPC algorithm is presented in Section \ref{STAMPC:sec:MPC_alg}, followed by the analysis of the theoretical properties. In Section \ref{STAMPC:sec:simulation}, a numerical example and comparison study are given. Finally, some concluding remarks are presented in Section \ref{STAMPC:sec:conclusion}.

\textit{Notation}: In this paper, we use the notations $\mathbb{R}, \mathbb{R}^{n}$ and $\mathbb{R}^{m\times n}$ to denote the sets of real numbers, column real vectors with $n$ elements and real matrices consisting of $n$ columns and $m$ rows, respectively. The set $\mathbb{B}^m=\{b\in\mathbb{R}^m: \|b\|_\infty\leq 1\}$ is called a unit hypercube of order $m$. Let $\mathbb{N}$ denote the set of non-negative integers, then sets $\mathbb{N}_{\geq a}$ and $\mathbb{N}_{[a,b]}$ represent  $\{x\in\mathbb{N}: x\geq a\}$ and  $\{x\in\mathbb{N}: a\leq x\leq b,b\geq a\}$, respectively. Similarly, the notation $\mathbb{R}_{\geq a}$ stands for the set $\{x\in\mathbb{R}: x\geq a\}$. Given a vector $x\in\mathbb{R}^{n}$, we use $\|x\|$ and $\|x\|_\infty$ to represent the Euclidean norm and infinity norm of $x$, respectively. Given two sets $X\subseteq \mathbb{R}^n$ and $Y\subseteq \mathbb{R}^n$, their Pontryagin difference is denoted by $X \ominus Y = \{z\in\mathbb{R}^n: z + y\in X; \forall y \in Y \}$, and their Minkowski sum is $X\oplus Y = \{x+y|x\in X, y\in Y \}$. A continuous and strictly increasing function $\alpha: \mathbb{R}_{\geq 0} \rightarrow \mathbb{R}_{\geq 0}$ is called a $\mathcal{K}$-function if $\alpha(0)=0$ and $\alpha(x)>0$ for all $x>0$. A function $\alpha: \mathbb{R}_{\geq 0} \rightarrow \mathbb{R}_{\geq 0}$ is called a $\mathcal{K}_\infty$-function if it is a $\mathcal{K}$-function and $\alpha(x)\rightarrow\infty$ as $x\rightarrow\infty$.
% % % % % % % % % % % % % % % % % % % % % % % % % % % % % % % % % % % % % % % % % % % % % % % % % % % % % % % % % % % % % % % % % % % % % % % % % % % % % %
\section{Preliminaries and problem formulation}\label{STAMPC:sec:problem}

Consider a discrete-time nonlinear system subject to parametric uncertainties and additive disturbances
\begin{equation}\label{STAMPC:sys_nonlinear_affine}
x_{t+1}=\mathcal{F}(x_t,u_t,v_t,d_t)\triangleq f(x_t,u_t)+g(x_t,u_t)v_t + d_t,
\end{equation} 
where $x_t\in\mathbb{R}^{n_x}, u_t\in\mathbb{R}^{n_u},v_t\in\mathbb{R}^{n_v}$ and $d_t\in\mathbb{R}^{n_x}$ are the system state, the control input, the time-varying parametric uncertainty and the additive disturbance, respectively. %Without loss of generality, $d_t$ and $v_t$ are possibly time-varying.
 $f:\mathbb{R}^{n_x}\times\mathbb{R}^{n_u}\rightarrow\mathbb{R}^{n_x}$ and $g:\mathbb{R}^{n_x}\times\mathbb{R}^{n_u}\rightarrow\mathbb{R}^{n_x\times n_{v}}$ are {known} nonlinear functions satisfying $f(0,0)=0$ and $g(0,0)=0$.  It is assumed that $x_t\in\mathcal{X},u_t\in\mathcal{U},v_t\in\mathcal{V}$ and $d_t\in\mathcal{D}$, where $\mathcal{X},~\mathcal{U}$ are compact sets and $\mathcal{V}$ and $\mathcal{D}$ are compact and convex polytopes. We also assume that  $\mathcal{X,V,D,U}$ contain the origin. 

Before presenting the main results, we firstly recall some well-established definitions  used in this paper.
 \begin{definition}[RPI set \cite{limon_input_2006Auto}]\label{STAMPC:def:RPI_set}
 	Consider a discrete-time uncertain system $x_{t+1} = \mathcal{G}(x_t,w_t)$, where $w_t\in\mathcal{W}$ denotes the model uncertainty and $x_t$ is the system state. A set $\Omega$ is a robust positively invariant (RPI) set for the system $x_{t+1} = \mathcal{G}(x_t,w_t)$ if $\mathcal{G}(x_t,w_t)\in\Omega$ for all $x_t\in\Omega$ and $w_t\in\mathcal{W}$. 
 \end{definition}
 
 \begin{definition}[Zonotope of order $n\times m$ \cite{wan_guaranteed_2018TAC}]\label{STAMPC:def:zonotope}
 	Given $p\in\mathbb{R}^n$ and $H\in\mathbb{R}^{n\times m}$, a zonotope of order $n\times m$ is a set of $n$-dimensional vectors defined by $\mathcal{Z} = p\oplus H\mathbb{B}^m = \{p+Hs:s\in\mathbb{B}^m\}.$
 \end{definition}
 
 \begin{definition}[ISpS-Lyapunov function \cite{limon_input_2006Auto}] \label{se2:def:ISpS_lyap}
 	%Suppose that Assumption \ref{sec2:assum:w_set} holds. 
 	Consider a nonlinear system in (\ref{STAMPC:sys_nonlinear_affine}).	A function $V(\cdot):\mathbb{R}^{n_x}\rightarrow\mathbb{R}_+$ is called ISpS-Lyapunov function if there exist $\mathcal{K}$-functions $\alpha_1(\cdot),\alpha_2(\cdot),\alpha_3(\cdot)$, constants $\gamma_1, \gamma_2$ and  a $\mathcal{K}$-function $\rho(\cdot)$ such that
 	\begin{subequations}
 		\begin{align}
 		&\alpha_1(\|x\|)\leq V(x)\leq \alpha_2(\|x\|)+\gamma_1,\label{STAMPC:ISPS_L_1}\\
 		&V(\mathcal{G}(x,v,d))-V(x)\leq -\alpha_3(\|x\|)+\rho(\|d\|)+\gamma_2.\label{STAMPC:ISPS_L_2}
 		\end{align}
 	\end{subequations}
 \end{definition}
 
 In the standard periodic MPC framework, a sequence of optimal control actions is obtained by solving the optimization problem at each time instant. Then the first element in this sequence will be sent to the actuator through the communication network and be implemented to the plant.  In order to reduce communication load, {a self-triggered adaptive  control is proposed: The unknown parameters are identified at sampling time instants. Then based on the new estimation of unknown parameters, the optimization problem is reformulated and solved to determine the control input and the triggering interval. Let $t_k$ denote the triggering time instant when the optimization problem needs to be solved. We consider the control policy and the set-based parameter estimator in the form of 
\begin{equation}\label{STAMPC:ut_gen}
\begin{array}{l}
u_t = \tau(x_{t_k},t-t_k),t\in\mathbb{N}_{[t_k,t_{k+1}-1]}\\
\hat{\mathcal{V}}_{t_{k+1}} = \Psi(\hat{\mathcal{V}}_{t_k},x_{t_k},x_{t_{k+1}},\mathbf{u}_{[t_k,t_{k+1}-1]},\mathcal{D})
\end{array}
\end{equation}
where $\hat{\mathcal{V}}_{t_k}$ is the EFSS of unknown parameters at time instant $t_k$ with  $\hat{\mathcal{V}}_{0} = \mathcal{V}$; $\mathbf{u}_{[t_k,t_{k+1}-1]} = \{ u_{t_k},u_{{t_k}+1},\cdots,u_{t_{k+1}-1}\}$. $t_k$ is obtained by using the following  self-triggering scheduler
\begin{equation}\label{STAMPC:st_scheduler}
t_{k+1} = t_k+H^*(x_{t_k},\hat{\mathcal{V}}_{t_k}),k\in\mathbb{N}_{>0},
\end{equation}
with $t_0 = 0$. Our objective is to design the control policy $\tau:\mathbb{R}^{n_x}\times\mathbb{N}\rightarrow\mathbb{R}^{n_u}$, the set-valued mapping $\Psi:\mathbb{R}^{n_v}\times\mathbb{R}^{n_x}\times\mathbb{R}^{n_x}\times\mathbb{R}^{(t_{k+1}-1-t_k)n_u}\times\mathbb{R}^{n_d}\rightrightarrows\mathbb{R}^{n_v}$ and the scheduling function $H^*:\mathbb{R}^{n_x}\times\mathbb{R}^{n_v}\rightarrow\mathbb{N}_{>0}$ such that: 1) The proposed set-membership parameter estimator provides a tight overestimation of the FSS; 2) the online parameter adaption is facilitated in the co-design of control inputs and triggering intervals to further reduce the average sampling frequency; 3) the system in (\ref{STAMPC:sys_nonlinear_affine}) is robustly stabilized with guaranteed robust constraint satisfaction and less conservative performance under recursive updates of the  system model.

% % % % % % % % % % % % % % % % % % % % % % % % % % % % % % % % % % % % % % % % % % % % % % % % % % % % % % % % % % % % % % % % % % % % % % % % % % % % % %
\section{Parameter estimation}\label{STAMPC:sec:estimation}
Define $y_t = x_{t}-f(x_{t-1},u_{t-1})$ and $\phi_{t-1} = g(x_{t-1},u_{t-1})$. According to (\ref{STAMPC:sys_nonlinear_affine}), we have the following regression model
	\begin{equation}\label{STAMPC:est_regression}
	y_t = \phi_{t-1} v_{t-1} + d_{t-1}.
	\end{equation}
Before presenting the parameter estimation algorithm, we firstly introduce the definition of the information set.
\begin{definition}[Information set]\label{STAMPC:def:infor_set}
	A set $\mathcal{L}_t$ is called the information set at time $t$ if it is consistent with the system in (\ref{STAMPC:est_regression}), the measurements $y_t$, the regressor $\phi_{t-1}$ and the uncertainty set $\mathcal{D}$, namely: $\mathcal{L}_t = \{{v_t}\in\mathbb{R}^{n_v}: y_t - \phi_{t-1} {v_t} \in\mathcal{D}\}.$
\end{definition}

Consider the dynamic evolution of model uncertainty $v_t$ with a general form  
\begin{equation}\label{STAMPC:sys_uncer_v}
v_{t+1} = \eta(v_t,\delta_t),
\end{equation}
where $\delta_t\in\mathcal{M}\subseteq\mathbb{R}^{n_\delta}$ is an auxiliary variable and $\eta:\mathbb{R}^{n_v}\times\mathbb{ R}^{n_\delta}\rightarrow\mathbb{R}^{n_v}$ is a nonlinear function. Let $\mathcal{V}_t$ denote the FSS of the uncertain parameter $v_t$ at time $t$. Given the information set $\mathcal{L}_{t+1}$, based on the system in (\ref{STAMPC:sys_uncer_v}), we have
\begin{equation}\label{STAMPC:set_exact_comp}
\mathcal{V}_{t+1} = \eta(\mathcal{V}_t,\mathcal{M})\cap\mathcal{L}_{t+1}.
\end{equation}
It can be seen from (\ref{STAMPC:set_exact_comp}) that the parameter estimation problem consists of two parts: 1) The set computation associated with the dynamic evolution, and 2) the intersection of the information set and the set obtained from the dynamic evolution. 
\begin{remark}
	The proposed method is derived based on the model of parametric uncertainty in (\ref{STAMPC:sys_uncer_v}). Its main benefit is to reduce the overestimation by effectively using the model information. But in many practical control problems, it may be difficult to find an exact model to describe $v_t$. Alternatively, we  consider the auxiliary variable $\delta_t$ in (\ref{STAMPC:sys_uncer_v}) to establish the uncertain model for $v_t$. For example, if $v_t$ in (\ref{STAMPC:sys_nonlinear_affine}) is subject to a bounded change rate $\bar{\delta}>0$, i.e., $\|v_{t+1}-v_t\|\leq\bar{\delta}$, we can find the following uncertain linear system  $v_{t+1} = \eta(v_t,\delta_t) \triangleq v_t + \delta_t$ with $\delta_t\in\mathcal{M}=\{\delta\in\mathbb{R}^{n_v}:\|\delta\|\leq\bar{\delta} \} $. In addition, the set $\mathcal{M}$ can also be time-varying, which can be handled directly via replacing $\mathcal{M}$ by $\mathcal{M}_t$ at sampling time instants.  
\end{remark}

As shown in (\ref{STAMPC:set_exact_comp}), the key issue of the parameter estimation problem is how to calculate the set-based dynamic evolution with less overestimation. In this section, we start by recalling the preliminary results on computing $\eta(\mathcal{V}_t,\mathcal{M})$ based on the indirect polytopic set computation \cite{wan_guaranteed_2018TAC}. Then a zonotope-based set-membership parameter estimator is presented for the nonlinear system in (\ref{STAMPC:sys_nonlinear_affine}) under the aperiodic sampling induced by the self-triggering scheme.
%If $\mathcal{V}_t$ is a zonotope, $\eta(\mathcal{V}_t,\mathcal{M})$ can be directly computed based on \textit{ K{\"u}hn's method} \cite{kuhn_rigorously_1998Computing,alamo_guaranteed_2005automatica,wan_guaranteed_2018TAC}. Hence, in this section, we start by recalling the preliminary results on the zonotopic set computation and the indirect polytopic set computation. Then according to this indirect method, a zonotope-based set-membership parameter estimator is presented for the nonlinear system in (\ref{STAMPC:sys_nonlinear_affine}) with aperiodic sampling.
% % % % % % % % % % % % % % % % % % % % % % % % % % % % % % % % % % % % % % %
\subsection{Indirect polytopic set computation}\label{STAMPC:sec:zono_computation}
%Although the set obtained from K{\"u}hn's method is a zonotope, intersecting between $\eta(\mathcal{V}_t,\mathcal{M})$ and $\mathcal{L}_{t+1}$ may result in a polytope set, which renders K{\"u}hn's method invalid in next time instant. To solve this problem, a novel indirect polytopic set computation method proposed in \cite{wan_guaranteed_2018TAC}. 
Let $\hat{\mathcal{V}}_t$ denote the EFSS for $v_t$ at time $t$. If $\hat{\mathcal{V}}_t$ is a zonotope with $\hat{\mathcal{V}}_t = p_t\oplus H_t \mathbb{B}^{n_{v_t}}$, the set-based dynamic evolution of the system in (\ref{STAMPC:est_regression}) can be computed by using the centered inclusion function \cite[Theorem 2]{alamo_guaranteed_2005automatica}, 
\begin{equation}\label{STAMPC:center_inclu_func} 
\eta_c(\hat{\mathcal{V}}_t,\mathcal{M}) = \eta(p_t,\mathcal{M})  + \triangledown_v\eta(\bar{\mathcal{V}}_t,\mathcal{M})(\hat{\mathcal{V}}_t\ominus p_t),
\end{equation}
where $\bar{\mathcal{V}}_t$ is a box bounding the set $\hat{\mathcal{V}}_t$. Assume that $\eta(p_t,\mathcal{M})$ is bounded by a zonotope $\check{p}_t\oplus\check{H}_t\mathbb{B}^{n_{\check{v}_t}}$, then we have
\begin{equation*}
\eta_c(\hat{\mathcal{V}}_t,\mathcal{M})\subseteq\check{p}_t\oplus\check{H}_t\mathbb{B}^{n_{\check{v}_t}}\oplus{\mathbb{M}}_t\mathbb{B}^{n_{v_t}}
\end{equation*} 
where $\hat{\mathbb{M}}_t = \triangledown_v\eta(\bar{\mathcal{V}}_t,\mathcal{M})H_t$. 
By using the zonotope inclusion operator shown in \cite[Theorem 3]{alamo_guaranteed_2005automatica}, we can further find a zonotope that bounds $\eta_c(\mathcal{V}_t,\mathcal{M})$.
\begin{remark}
	As proposed in \cite{alamo_guaranteed_2005automatica}, the zonotope bounding the function  $\eta(p_t,\mathcal{M})$ can be found by using a natural interval extension or a mean value extension. Furthermore, if $\mathcal{M}$ is a zonotope with $\mathcal{M} = p_{\delta}\oplus H_\delta\mathbb{B}^{n_\delta}$, and  $\eta(v_t,\delta_t)$ is an affine function of the disturbance $\delta_t$, i.e., $\eta(v_t,\delta_t) =  \eta_1(v_t)+\eta_2(v_t)\delta_t$, we can find a zonotope $p_{as}\oplus H_{as}\mathbb{B}^{n_\delta}$ such that $\eta(v_t,\mathcal{M}) = p_{as}\oplus H_{as}\mathbb{B}^{n_\delta}$ where $p_{as} =\eta_1(p_t)+\eta_2(p_t)p_\delta$ and $H_{as} = \eta_2(p_t)H_{\delta}$.
\end{remark}

Although the set obtained from the center inclusion function (\ref{STAMPC:center_inclu_func}) is a zonotope, performing the intersection operation in (\ref{STAMPC:set_exact_comp}) may result in a polytope, rendering the center inclusion function (\ref{STAMPC:center_inclu_func}) infeasible at next time instant. Hence, in the conventional zonotope-based set-membership methods, e.g., \cite{alamo_guaranteed_2005automatica,wang_zonotope_2018Automatica}, a zonotopic bounding process is imposed on the set resulted from the intersection. But this bounding process may lead to unnecessary overestimation. To avoid this, we recall the novel indirect polytopic set computation technique \cite{wan_guaranteed_2018TAC} in the following.

Suppose that $\hat{\mathcal{V}}$ is a polytope and $\hat{\mathcal{V}} = \mathcal{Z}_1\cap\mathcal{Z}_2$, where $\mathcal{Z}_1$ and $\mathcal{Z}_2$ are two polytopes. According to set theory, it can be derived that $\eta(\hat{\mathcal{V}}, \mathcal{M}) = \eta(\mathcal{Z}_1\cap\mathcal{Z}_2, \mathcal{M})\subseteq \eta(\mathcal{Z}_1, \mathcal{M})\cap\eta(\mathcal{Z}_2, \mathcal{M})$. If $\mathcal{Z}_1$ and $\mathcal{Z}_2$ are zonotopic,
$\eta(\mathcal{Z}_1, \mathcal{M})$ and $\eta(\mathcal{Z}_2, \mathcal{M})$ can be computed by using the center inclusion function (\ref{STAMPC:center_inclu_func}) directly without the bounding process. Therefore, the essence of indirect polytopic set computation is to find a group of zonotopes whose intersection is the polytope. It is worth mentioning that the results in \cite{wan_guaranteed_2018TAC} consider a 2-dimensional case; the following lemma extends this useful technique to accommodate general $n$-dimensional cases.}
\begin{lemma}
	Let $\mathcal{X}\in\mathbb{R}^{n}$ denote a convex and compact polytope. Suppose that $\mathcal{X}$ consists of $m$ half-spaces. Then the polytope $\mathcal{X}$ can be exactly represented by the intersection of at most $m$ zonotopic sets.
\end{lemma}

\begin{IEEEproof}
	Let $h_i,i\in\mathbb{N}_{[1,m]}$ denote the $i$th half-space in $\mathcal{X}$. By using $h_i$ as an edge, a parallelotope or a zonotope $\mathcal{H}_i$ can be constructed such that $\mathcal{X}\subseteq \mathcal{H}_i$. Therefore, we can find a sequence of zonotopes $\{\mathcal{H}_i\}_{i\in\mathbb{N}_{[1,m]}}$ whose intersection is $\mathcal{X}$. In addition, if the polytope is symmetric, there exist $i,j\in\mathbb{N}_{[1,m]}$ such that $\mathcal{X}\subseteq\mathcal{H}_i\subseteq \mathcal{H}_j$. Consequently, at most $m$ zonotopes are needed to construct the polytope.
\end{IEEEproof}

\begin{remark}
	For each $h_i$, we can find a point $\bar{x}_i\in\mathcal{X}$ by solving a linear programming problem such that distance between $x_i$ and $h_i$ is maximized. According to
	$x_i$, we can find a half-space $\bar{h}_i$ that is parallel to $h_i$. 	Suppose that $\mathcal{X}\triangleq\texttt{conv}(\{\hat{x}_l\}_{l\in\mathbb{N}_{[0,\bar{n}]}})$, where $\texttt{conv}(\cdot)$ is the convex hull of a set of points; $\hat{x}_l$ denotes the $l$th vertex of $\mathcal{X}$; $\bar{n}$ is the amount of vertexes in $\mathcal{X}$. Let
	$\mathcal{N}_i = \{l|i\in\mathbb{N}_{[1,m]},l\in\mathbb{N}_{[0,\bar{n}]}  \}$ denote the set of indices of vertex located at $h_i$, and $\tilde{x}_l$ is the projection of $\hat{x}_l$ onto $\bar{h}_i$ with $l\in\mathcal{N}_i$. Then the zonotope $\mathcal{H}_i$ can be constructed as $\mathcal{H}_i=\texttt{conv}(\{\hat{x}_l,\tilde{x}_l\}_{l\in\mathcal{N}_i})$. 
	%The construction of zonotopes for 2-dimensional polytopes is investigated in \cite{wan_guaranteed_2018TAC}, where it is proved that the number of zonotopes needed to form the polytope is $\frac{m}{2}$ if $m$ is even or $\frac{m+1}{2}$ if  $m$ is odd.
\end{remark}

\subsection{Set-membership parameter estimation under the self-triggering mechanism}
Without loss of generality, it is assumed that the triggering interval computed at time $t_k$ is $H^*_{t_k}$. If $H^*_{t_k} = 1$, according to Definition \ref{STAMPC:def:infor_set}, the information set $\mathcal{L}_{t_k}$ can be calculated as follows
\begin{equation}\label{STAMPC:infoSet_consecutive}
\begin{array}{l}
\mathcal{L}_{t_k} = \{{v_{t_k}}\in\mathbb{R}^{n_v}: x_{t_k}-f(x_{t_k-1},u_{t_k-1}) \\
~~~~~~~~~~~~~~~~~~~~~~~- g(x_{t_k-1},u_{t_k-1}) {v_{t_k}} \in\mathcal{D}\}.
\end{array}
\end{equation}
However, under the self-triggering scheduler (\ref{STAMPC:st_scheduler}), the system states at two consecutive time instants become inaccessible if $H^*_{t_k} > 1$, making the calculation of $\mathcal{L}_{t_k}$ in (\ref{STAMPC:infoSet_consecutive}) infeasible. To solve this problem, we propose to 1) estimate the reachable sets of system states between two sampling instants, and 2) approximate the information sets at triggering time instants by using the estimated reachable set instead of the measured state.

Suppose that $H^*_{t_k} > 1$. Since all of the system states are accessible at triggering time instants, we define $X_{t_k} \triangleq \{x_{t_k}\} = p^x_{t_k}\oplus O^x_{t_k}\mathbb{B}^{n^x_{t_k}}$ with $p^x_{t_k} = x_{t_k}, n^x_{t_k} = 0$ and $O^x_{t_k} = \mathbf{0}$. Let  $\hat{{X}}_{l+t_k}\triangleq p^x_{l+t_k}\oplus O^x_{l+t_k}\mathbb{B}^{n^x_{l+t_k}}$ denote the estimated reachable set for the system state $x_{l+t_k}$ from $X_{t_k}$, where $l\in\mathbb{N}_{[0,H^*_{t_k}-1]}$ and $\hat{X}_{t_k} = X_{t_k}$.
Then by using the zonotopic set computation method in Section \ref{STAMPC:sec:zono_computation}, $\hat{X}_{l+1+t_k}$ can be estimated as follows:
\begin{equation}\label{STAMPC:open_state_est}
\hat{{X}}_{l+1+t_k} = \mathcal{F}_c(\hat{{X}}_{l+t_k},u_{l+t_k},\hat{\mathcal{V}}_{l+t_k},\mathcal{D})
\end{equation} 
where the function $\mathcal{F}_c(\hat{{X}}_{l+t_k},u_{l+t_k},\hat{\mathcal{V}}_{l+t_k},\mathcal{D}) \triangleq  \triangledown_x\mathcal{F}(\bar{{X}}_{l+t_k},u_{l+t_k},\hat{\mathcal{V}}_{l+t_k},\mathcal{D})(\hat{{X}}_{l+t_k} - p^x_{l+t_k})+ \mathcal{F}(p^x_{l+t_k},u^*_{l+t_k},\hat{\mathcal{V}}_{l+t_k},\mathcal{D})\nonumber$
is the centered inclusion function derived in \cite[Theorem 3]{alamo_guaranteed_2005automatica}; $\bar{{X}}_{l+t_k}$ is a box bounding $\hat{{X}}_{l+t_k}$. Since $x_{t_k+1}$ is inaccessible, we consider the relatively conservative estimation of $\mathcal{V}_{t_k+l}$ in (\ref{STAMPC:open_state_est}) by ignoring the intersecting operation
\begin{equation}\label{STAMPC:over_est_V}
\hat{\mathcal{V}}_{t_k+l+1} = (\cap_{i=0}^{n_{l}}\eta_c(\check{\mathcal{V}}_{i,l},\mathcal{M}))\cap \mathcal{V}
\end{equation}
where $i\in\mathbb{N}_{[1,n_{t_k}]},~l\in\mathbb{N}_{[0,H^*_{t_k}-2]}$, and $\{\check{\mathcal{V}}_{i,l+t_{k}}\}$ is the sequence of zonotopes satisfying $\hat{\mathcal{V}}_{l+t_k} = \cap\{\check{\mathcal{V}}_{i,l+t_{k}}\}$. Consequently, based on the system in (\ref{STAMPC:sys_nonlinear_affine}), the information set at time $t_{k+1}$ can be approximated by 
\begin{align}\label{STAMPC:est_info_set}
\hat{\mathcal{L}}_{t_{k+1}} = \{&v_{t_{k+1}}\in\mathbb{R}^{n_v}: x_{t_{k+1}}-f(x,u_{t_{k+1}-1}) \nonumber\\ & -g(x,u_{t_{k+1}-1})v_{t_{k+1}}\in\mathcal{D},x\in\hat{{X}}_{t_{k+1}-1}\}.
\end{align} 
Therefore, we have
\begin{equation}\label{STAMPC:over_est_V_final}
\hat{\mathcal{V}}_{t_{k+1}} =  (\cap_{i=0}^{n_{l}}\eta_c(\check{\mathcal{V}}_{i,{t_{k+1}-1}},\mathcal{M}))\cap\hat{\mathcal{L}}_{t_{k+1}}\cap \mathcal{V}.
\end{equation}
In summary, the procedure for updating the EFSS $\hat{\mathcal{V}}_{t+1}$ is presented in Algorithm \ref{STAMPC:alg:set_update}.

\begin{algorithm}[!t]
	\caption{Zonotope-based set-membership parameter estimation algorithm}
	\begin{algorithmic}[1]\label{STAMPC:alg:set_update}
		\REQUIRE Measured system states $x_{t_{k-1}}$ and $x_{t_{k}}$; sampling time instants $t_{k-1}$ and $t_k$; control input sequence $\{u_{i}\},i\in\mathbb{N}_{[t_{k-1},t_k-1]}$; EFSS $\hat{\mathcal{V}}_{t_{k-1}}$; uncertainty set $\mathcal{D}$. 
		\STATE Find a sequence of zonotopes $\{\check{\mathcal{V}}_{i,t_{k-1}}\}$ such that $\hat{\mathcal{V}}_{t_{k-1}} = \cap_{i=1}^{n_{t_{k-1}}}
		\{\check{\mathcal{V}}_{i,t_{k-1}}\}$; Set $\hat{X}_{t_{k-1}} = \{x_{t_{k-1}}\}$.\;
		\IF{$t_{k}-t_{k-1} > 1$}
		\FOR {$j=t_{k-1},t_{k-1}+1,\cdots,t_{k}-2$}
		% \STATE Calculate $\eta_c(\check{\mathcal{V}}_{i,j},\mathcal{M})$ by using (\ref{STAMPC:center_inclu_func}), and then get $\hat{\mathcal{V}}_{j+1} = \cap\eta_c(\check{\mathcal{V}}_{i,j},\mathcal{M})$
		\STATE Estimate $\hat{{X}}_{j+1}$ and $\hat{\mathcal{V}}_{j+1}$ by using (\ref{STAMPC:open_state_est}) and (\ref{STAMPC:over_est_V}), respectively.
		\ENDFOR
		\STATE Compute ${\mathcal{L}}_{t_k}$ by using (\ref{STAMPC:est_info_set}), then $\hat{\mathcal{V}}_{t_{k}}$ can be obtained by following  (\ref{STAMPC:over_est_V_final}).
		\ELSE
		\STATE Calculate the information set $\mathcal{L}_{t_k}$ by using (\ref{STAMPC:infoSet_consecutive}), and then compute the new EFSS $\hat{\mathcal{V}}_{t_k} = \eta_c(\check{\mathcal{V}}_{1,t_{k-1}},\mathcal{M})\cap\eta_c(\check{\mathcal{V}}_{2,t_{k-1}},\mathcal{M})\cap\cdots\eta_c(\check{\mathcal{V}}_{n_{t_{k-1}},t_{k-1}},\mathcal{M})\cap\mathcal{L}_{t_{k}}$. 	
		\ENDIF 
	\end{algorithmic}
\end{algorithm}

% % % % % % % % % % % % % % % % % % % % % % % % % % % % % % % % % % % % % % % % % % % % % % % % % % % % % % % % % % % % % % % % % % % % % % % %
\section{Self-triggered adaptive min-max MPC}\label{STAMPC:sec:MPC_alg}
{In this section, based on the set-membership parameter estimation scheme described in Section \ref{STAMPC:sec:estimation}, we firstly present the min-max MPC optimization problem. Thereafter, a self-triggering scheduler accommodating the estimation of uncertainty is proposed, followed by a summary of the proposed self-triggered adaptive min-max MPC algorithm. Finally, this section concludes with an analysis of closed-loop stability and recursive feasibility.}

\subsection{Min-max optimization}
Given the prediction horizon $N\in\mathbb{N}_{\geq 0}$ and an integer $H\in\mathbb{N}_{[1,N]}$, we formulate the following cost function for the MPC problem at time $t_k$
\begin{equation*}
\begin{array}{l}
J_N^H(x_{t_k},\mathbf{u}_{t_k,N},\mathbf{v}_{t_k,N},\mathbf{d}_{t_k,N},\beta_{t_k}) \triangleq \sum_{l=0}^{H-1}\frac{1}{\beta_{t_k}}\ell(x_{l|t_k},u_{l|t_k})\\
+\sum_{l=H}^{N-1}\ell(x_{l|t_k},u_{l|t_k})+ \ell_f(x_{N|t_k}),
\end{array}
\end{equation*}
where $\beta_{t_k}\in\mathbb{R}_{\geq 1}$ is a scalar to be designed, $x_{l|t_k}$ is the predicted system state $l$ steps ahead from the time $t_k$ with $l\in\mathbb{N}_{[0,N]}$. $\mathbf{u}_{t_k,N}$ denotes the control sequence and  $\mathbf{u}_{t_k,N}=\{u_{0|t_k},u_{1|t_k},\cdots,u_{N-1|t_k}\}$. $\mathbf{v}_{t_k,N}=\{v_{1|t_k},v_{2|t_k},\cdots,v_{N|t_k}\}$ and $\mathbf{d}_{t_k,N}=\{d_{1|t_k},d_{2|t_k},\cdots,d_{N|t_k}\}$ are sequences of multiplicative and additive disturbances, respectively. We suppose that the stage cost function $\ell:\mathbb{R}^{n_x}\times\mathbb{R}^{n_u}\rightarrow\mathbb{R}_{\geq0}$ and the terminal cost function $\ell_f:\mathbb{R}^{n_x}\rightarrow\mathbb{R}_{\geq0}$ are continuous functions satisfying $\ell(0,0)=0$ and $\ell_f(0)=0$. 

In order to guarantee robust constraint satisfaction, we propose a self-triggered adaptive min-max MPC approach. Let $\mathbf{u}_{t_k,N}$ denote the decision variable for the MPC optimization problem, then the proposed adaptive MPC algorithm is based on the following min-max optimization problem
\begin{subequations}\label{STAMPC:opt}
	\begin{align}
	&V_N^H(x_{t_k},\beta_{t_k})=
	\min_{u_{l|t_k}\in\mathcal{U},l\in\mathbb{N}_{[0,H-1]}}\{\max_{v_{l|t_k}\in\hat{\mathcal{V}}_{l|t_k}, d_{l|t_k}\in\mathcal{D}}\nonumber\\
	&\{\sum_{l=0}^{H-1}\frac{1}{\beta_{t_k}}\ell(x_{l|t_k},u_{l|t_k})+V_{N-H}(x_{H|t_k})\}~\text{such that} \nonumber\\
	&x_{H|t_k}\in\mathcal{X}_{N-H},\forall~v_{l|t_k}\in\hat{\mathcal{V}}_{l|t_k}, d_{l|t_k}\in\mathcal{D},l\in\mathbb{N}_{[0,H-1]}    \},\\
	&\text{s.t.}~x_{l+1|t_k} = \mathcal{F}(x_{l|t_k},u_{l|t_k},v_{l|t_k},d_{l|t_k}),l\in\mathbb{N}_{[0,H-1]},\\
	&~~~~x_{H|t_k}\in\mathcal{X}_{N-H},\\
	&~~~~x_{0|t_k} = x_{t_k},x_{l|t_k}\in\mathcal{X},l\in\mathbb{N}_{[0,H-1]},
	\end{align}
\end{subequations}
where
\begin{align}\label{STAMPC:opt_closed}
&V_{N-i}(x_{i|t_k}) = \min_{u_{i|t_k}\in\mathcal{U}}\{\max_{v_{i|t_k}\in\hat{\mathcal{V}}_{i|t_k}, d_{i|t_k}\in\mathcal{D}} \{\ell(x_{i|t_k},u_{i|t_k})\nonumber\\
&+V_{N-i-1}(\mathcal{F}(x_{i|t_k},u_{i|t_k},v_{i|t_k},d_{i|t_k}))\}~\text{such that}\nonumber\\ &\mathcal{F}(x_{i|t_k},u_{i|t_k},v_{i|t_k},d_{i|t_k})\in\mathcal{X}_{N-i-1},\forall v_{i|t_k}\in\hat{\mathcal{V}}_{i|t_k},\nonumber\\
&\forall d_{i|t_k}\in\mathcal{D}\},
\end{align}
$i\in\mathbb{N}_{[H,N-1]}$ and $\mathcal{X}_{N-i}$ is the set of admissible states which can be robustly steered into the terminal set $\mathcal{X}_f$ within $N-i$ steps. $\hat{\mathcal{V}}_{i|t_k}$ is the predicted EFSS $i$ step ahead from the time $t_k$ with $\hat{\mathcal{V}}_{0|t_k} = \hat{\mathcal{V}}_{t_k}$. In addition, since $v_t\in\mathcal{V}$ for all $t\geq 0$, $\hat{\mathcal{V}}_{i|t_k}$ can be calculated recursively by following
\begin{equation}\label{STAMPC:v_cen_incl_func}
\hat{\mathcal{V}}_{i+1|t_k}=
(\cap_{s=0}^{n_{l}}\eta_c(\check{\mathcal{V}}_{s,l|t_k},\mathcal{M}))\cap \mathcal{V}, ~~~~i\in \mathbb{N}_{[0,N-1]}
\end{equation}
where $\hat{\mathcal{V}}_{i|t_k} = \cap_{s=0}^{n_{l}} \{\check{\mathcal{V}}_{s,i|t_k}\}$. Furthermore, the initial conditions are specified as $V_0(x_{N|t_k}) \triangleq \ell_f(x_{N|t_k})$ and $\mathcal{X}_0\triangleq\mathcal{X}_f$, respectively.

{As shown in the min-max optimization problem (\ref{STAMPC:opt}), $\beta_{t_k}$ determines the penalty on the cost for open loop scenarios. Hence, it is possible to achieve a less conservative closed-loop performance by designing suitable $\beta_{t_k}$. In order to further reduce the average sampling frequency, we propose the evolution of $\beta_{t_k}$ by implicitly considering the EFSS $\hat{\mathcal{V}}_{t_k}$ in the following.
% % % % % % % % % % % % % % % % % % % % % % % % % % % % % % % % % % % % % % %
\subsection{Self-triggering scheduler} 
Under the self-triggering mechanism,
the MPC problem is solved at the triggering time instant only. Therefore, between two consecutive triggering time instants, the system in (\ref{STAMPC:sys_nonlinear_affine}) is operated with the open loop control actions determined by (\ref{STAMPC:opt}), i.e., $u_t = \tau(x_{t_k},t-t_k) = u^*_{t-t_k|t_k},t\in\mathbb{N}_ {[t_k,t_{k+1}-1]}$ and $\{u^*_{t-t_k|t_k}\}$ denotes the optimal solution of the MPC optimization problem (\ref{STAMPC:opt}) obtained at time $t_k$. Inspired by \cite{acipl_liu_robust_2019Auto}, we design the following self-triggering scheduler
\begin{equation}\label{STAMPC:opt:self_trig_cond}
\begin{array}{l}
t_{k+1} = t_k + H^*(x_{t_k},\hat{\mathcal{V}}_{t_k})\\
H^*(x_{t_k},\hat{\mathcal{V}}_{t_k}) \triangleq\max\{H\in\mathbb{N}_{[1,H_{\max}]}|\\
~~~~~~~~~~~~~~V_N^H(x_{t_k},\beta_{t_k})\leq V_N^1(x_{t_k},\beta_{t_k})\}
\end{array}
\end{equation}
where $H_{\max}\in\mathbb{N}_{[1, {N}]}$ denotes the maximum number of time instants allowed for the open loop scenario. Consequently, the system in (\ref{STAMPC:sys_nonlinear_affine}) becomes
\begin{equation}\label{STAMPC:sys_closed_loop}
x_{t+1}=\mathcal{F}(x_t,\tau(x_{t_k},t-t_k),v_t,d_t).
\end{equation}

As shown in (\ref{STAMPC:opt:self_trig_cond}), the performance of the proposed self-triggered adaptive MPC algorithm depends on the parameter $\beta_{t_k}$. Using a larger $\beta_{t_k}$ will result in a larger triggering interval $H^*(x_{t_k},\hat{\mathcal{V}}_{t_k})$, but will deteriorate the regulation performance \cite{acipl_liu_robust_2019Auto}. Furthermore, the triggering interval is also determined by the bound of uncertainties: A larger bound of uncertainties renders a shorter triggering interval by following (\ref{STAMPC:opt:self_trig_cond}). To further reduce the average sampling frequency, we propose the following adaptive mechanism associated with $\hat{\mathcal{V}}_{t_k}$
\begin{equation}\label{STAMPC:beta_k}
\beta_{t_k} = \min(\frac{\xi_0}{\xi_{t_k}}\beta_0, \beta_{\max}),
\end{equation} 
where $\xi_{t_k} = \sup_{v_1,v_2\in \hat{\mathcal{V}}_{t_k}}\|v_1-v_2\|$; $\beta_0\in\mathbb{R}_{\geq 1}$ is a fixed constant, and $\beta_{\max}$ is a positive number. Since $H^*(x_{t_k},\hat{\mathcal{V}}_{t_k})$ is bounded by the prediction horizon $N$, increasing $\beta_{t_k}$ cannot enlarge the triggering intervals if $\beta_{t_k}$ is sufficiently large. Hence we manually bound $\beta_k$ by $\beta_{\max}$ for all $t_k\geq 0$. 

\begin{remark}
	To design a self-triggered adaptive MPC algorithm, the offline constructed RPI set $\mathcal{X}_f$ is commonly used to guarantee robust stability in the presence of additive and multiplicative uncertainties. At each triggering time instant, it is possible to improve the closed-loop behaviors by recomputing the terminal set $\mathcal{X}_f$ in (\ref{STAMPC:opt}) based on the new estimated uncertainty set. However, as remarked in \cite{lorenzen_robust_2019Automatica}, if the local control policy $\kappa_f$ were not simultaneously updated based on the new EFSS, the closed-loop performance would not be significantly improved. On the other hand, updating both $\mathcal{X}_f$ and the feedback gain would render the MPC problem much more complicated. Consequently, we consider the fixed terminal set $\mathcal{X}_f$ in our method, which is designed offline based on the initial uncertainty sets $\mathcal{V}$ and $\mathcal{D}$.
\end{remark}

According to the developed set-membership parameter estimator, the proposed self-triggered adaptive MPC scheme is summarized in Algorithm \ref{STAMPC:alg:stampc_state_para_est}.
\begin{algorithm}[!t]
	\caption{Self-triggered adaptive MPC with state and parameter estimation}
	\begin{algorithmic}[1]\label{STAMPC:alg:stampc_state_para_est}
		\REQUIRE Initial system state $x_{0}$; initial EFSS $\hat{\mathcal{V}}_{0}$; uncertainty set $\mathcal{D}$; tuning parameters $\beta_0$ and $H_{\max}$.
		\STATE Set $t =0,t_k = 0$ and $k = 0$. 
		\WHILE{The control action is not stopped}
		\STATE Measure the state $x_{t_k}$ of the system in (\ref{STAMPC:sys_nonlinear_affine}).
		\IF{$t_k>0$}
		\STATE Compute the EFSS $\hat{\mathcal{V}}_{t_k}$ by following \textbf{Algorithm \ref{STAMPC:alg:set_update}}.
		\ENDIF
		\STATE Reformulate and solve the optimization problems in (\ref{STAMPC:opt}) based on $\hat{\mathcal{V}}_{t_k}$ to obtain the optimal control sequence $\mathbf{u}_{t_k,N}^*$ and determine the next sampling time instant $t_{k+1} = t_k+H^*(t_k)$ by using (\ref{STAMPC:opt}) and (\ref{STAMPC:opt:self_trig_cond}).
		\STATE Implement the control input $u_{t_k+i} = u^*_{i|t_k}$ to the system in (\ref{STAMPC:sys_nonlinear_affine}) where $i\in\mathbb{N}_{[0,H^*(x_{t_k},\hat{\mathcal{V}}_{t_k})-1]}$.
		\STATE Set $t = t+H^*(x_{t_k},\hat{\mathcal{V}}_{t_k}), k = k +1$ and $t_k = t$.
		\ENDWHILE		
	\end{algorithmic}
\end{algorithm}
% % % % % % % % % % % % % % % % % % % % % % % % % % % % % % % % % % % % % % % % % % % % % % % % % % % % % % % % % % % % % % % % % % % % % % % % % % % % % %
\subsection{Recursive feasibility and closed-loop stability}
To ensure the closed-loop stability and recursive feasibility, we make the following standard assumptions.
\begin{assumption}\label{STAMPC:assump:gen_asum}
	There exist a local controller $\kappa_f:\mathbb{R}^{n_x}\rightarrow\mathbb{R}^{n_u}$, a $\mathcal{K}$-function $\alpha$ and some constants $\sigma_1,\sigma_2$ and $\sigma_3$ such that 
	\begin{itemize}
		\item[1)] $\mathcal{X}_f\subseteq\mathcal{X}$ is an RPI set containing the origin for the closed-loop system $x_{t+1}=\mathcal{F}(x_t,\kappa_f(x_t),v_t,d_t)$.
		\item[2)] $\kappa_f(x_t)\subseteq\mathcal{U}$ for all $x\in\mathcal{X}_f$.
		\item[3)] $\ell(x_t,u_t)\geq \sigma_1\|x_t\|$ for all $x_t\in\mathcal{X}$ and $u_t\in\mathcal{U}$.
		\item[4)] $\sigma_2\|x_t\|\leq \ell_f(x_t)\leq \sigma_3\|x_t\|$ for all $ x_t\in\mathcal{X}_f$.
		\item[5)] $\ell_f(\mathcal{F}(x_t,\kappa_f(x_t),v_t,d_t))-\ell_f(x_t)\leq -\ell(x_t,\kappa_f(x_t))+\alpha(\|d_t\|)$ for all $x_t\in\mathcal{X}_f,v_t\in\mathcal{V}$ and $d_t\in\mathcal{D}$.
	\end{itemize}
\end{assumption}

\begin{theorem}\label{STAMPC:thm:feasibility}
	For the nonlinear system in (\ref{STAMPC:sys_nonlinear_affine}), suppose that Assumption \ref{STAMPC:assump:gen_asum} holds, then the proposed self-triggered adaptive MPC scheme, as presented in Algorithm \ref{STAMPC:alg:stampc_state_para_est} is recursively feasible for all $x_0\in\mathcal{X}_N$. 
\end{theorem}

\begin{IEEEproof}
	Let $\hat{\mathcal{V}}_{t_k}$ denote the EFSSs obtained by Algorithm \ref{STAMPC:alg:stampc_state_para_est}. At time $t_k$, suppose that $x_{t_k}\in\mathcal{X}_N$, and the next sampling time instant is $t_{k+1} = t_k + H^{*}(x_{t_k},\hat{\mathcal{{V}}}_{t_k})$. For simplicity, we use the notation $H_{t_k}^*$ to denote $H^{*}(x_{t_k},\hat{\mathcal{{V}}}_{t_k})$ in the following. Let  $\mathbf{u}^*_{t_k,N}=\{u^*_{0|t_k},u^*_{1|t_k},\cdots,u^*_{N-1|t_k}\}$ denote the optimal solution of the optimization problem (\ref{STAMPC:opt}) at time $t_k$, then we construct the candidate input sequence at time $t_{k+1}:\bar{\mathbf{u}}_{t_{k+1},N} = \{u^*_{H^*(t_k)|t_k},u^*_{H^*(t_k)+1|t_k},\cdots,u^*_{N-1|t_k},\kappa_f(x_{N|t_k}),	\cdots,\\ \kappa_f(x_{N+H^*(t_k)-1|t_k})\}.$  
%	\begin{align*}
%	&\bar{\mathbf{u}}_{t_{k+1},N} = \{u^*_{H^*(t_k)|t_k},u^*_{H^*(t_k)+1|t_k},\cdots,u^*_{N-1|t_k},\kappa_f(x_{N|t_k}),\\
%	&\cdots,\kappa_f(x_{N+H^*(t_k)-1|t_k})\}.
%	\end{align*}
	 It follows from (\ref{STAMPC:opt_closed}) that each element in $\bar{\mathbf{u}}_{t_{k+1},N}$ is a feedback control policy. Since $\mathcal{X}_f$ is an RPI set for the system in (\ref{STAMPC:sys_nonlinear_affine}) and $\hat{\mathcal{V}}_{i|t_k}\subseteq\mathcal{V}$ for all $t_k\in\mathbb{N}_{\geq 0}$ and $i\in\mathbb{N}_{[0,N]}$, $\mathcal{X}_N$ is an invariant set for the system in (\ref{STAMPC:sys_nonlinear_affine}) under the proposed self-triggered adaptive MPC algorithm for all $v_{t}\in\mathcal{V}_t$ and  $d_{t}\in\mathcal{D}$. Hence the recursive feasibility is proved.
\end{IEEEproof}

\begin{remark}
As shown in Theorem \ref{STAMPC:thm:feasibility}, we investigate the recursive feasibility of the proposed method at triggering time instants. Although the control inputs are executed in an open loop configuration between two triggering time instants, it still can be guaranteed that $x_t\in\mathcal{X}$ and $u_t\in\mathcal{U}$ for all $t\geq 0$ and $x_0\in\mathcal{X}_N$ by considering the  worst-case realization of uncertainty in the MPC optimization problem.
\end{remark}
\begin{theorem}\label{STAMPC:thm:stability}
	For the nonlinear system in (\ref{STAMPC:sys_nonlinear_affine}), suppose that Assumption \ref{STAMPC:assump:gen_asum} holds and $x_0\in\mathcal{X}_N$. If the triggering time instants $t_k$ are determined by (\ref{STAMPC:opt:self_trig_cond}), then the closed-loop system in (\ref{STAMPC:sys_closed_loop}) is ISpS at the triggering time instants. 
\end{theorem}
\begin{IEEEproof}
	To prove ISpS of the closed-loop system, we need to show that the optimal cost function is an ISpS Lyapunov function. Next we will demonstrate the satisfaction of (\ref{STAMPC:ISPS_L_1}) and (\ref{STAMPC:ISPS_L_2}) for the optimal cost function $V_N^{H_{t_k}^*}(x_{t_k},\beta_{t_k})$.
	
	By applying Algorithm \ref{STAMPC:alg:stampc_state_para_est}, it follows from (\ref{STAMPC:opt}) that
	\begin{align*}
	V_N^{H_{t_k}^*}(x_{t_k},\beta_{t_k}) =& J_N^H(x_{t_k},\mathbf{u}^*_{t_k,N},\mathbf{v}^*_{t_k,N},\mathbf{d}^*_{t_k,N},\beta_{t_k})\\
	&\geq J_N^H(x_{t_k},\mathbf{u}^*_{t_k,N},\mathbf{0},\mathbf{0},\beta_{t_k})\\
	&\geq \min_{\mathbf{u}^*_{t_k}}J_N^H(x_{t_k},\mathbf{u}_{t_k,N},\mathbf{0},\mathbf{0},\beta_{t_k})\\
	%& \geq \frac{\sigma_1}{\beta_{t_k}}\|x_{t_k}\| \\
	&\geq \frac{\sigma_1}{\beta_{\max}}\|x_{t_k}\|.
	\end{align*}
	Hence the lower bound of $	V_N^{H_{t_k}^*}(x_{t_k},\beta_{t_k})$ is derived. Furthermore, let $\check{\mathbf{u}}_{t_k,N}$ denote the optimal solution associated with $V_N^{1}(x_{t_k},\beta_{t_k})$, and define a control sequence $\tilde{\mathbf{u}}_{t_k,N+1}$ as $\tilde{\mathbf{u}}_{t_k,N+1} = \{ \check{\mathbf{u}}_{t_k,N},\kappa_f(x_{N|t_k}) \} $. As such,  we have 
	\begin{align*}
	&J_{N+1}^1(x_{t_k},\tilde{\mathbf{u}}_{t_k,N+1},\mathbf{v}_{t_k,N+1},\mathbf{d}_{t_k,N+1},\beta_{t_k}) =\\ &~~~~~~~~J_{N}^1(x_{t_k},\check{\mathbf{u}}_{t_k,N},\mathbf{v}_{t_k,N},\mathbf{d}_{t_k,N},\beta_{t_k}) - \ell_f(x_{N|t_k})\\
	&~~~~~~~~+ \ell_f(x_{N+1|t_k}) + \ell(x_{N|t_k},\kappa_f(x_{N|t_k}). 
	\end{align*}
	for all $x_{t_k}\in\mathcal{X}_N$. Then it follows from Assumption \ref{STAMPC:assump:gen_asum} that 
	\begin{align*}
	&J_{N+1}^1(x_{t_k},\tilde{\mathbf{u}}_{t_k,N+1},\mathbf{v}_{t_k,N+1},\mathbf{d}_{t_k,N+1},\beta_{t_k})\\
	&~~~~\leq J_{N}^1(x_{t_k},\check{\mathbf{u}}_{t_k,N},\mathbf{v}_{t_k,N},\mathbf{d}_{t_k,N},\beta_{t_k})+\alpha(\|d_{t_k}\|).
	\end{align*}
	According to the suboptimality of the sequence $\tilde{\mathbf{u}}_{t_k,N+1}$, it is inferred that 
	\begin{align}\label{STAMPC:V_N_Mono}
	&V_{N+1}^1(x_{t_k},\beta_{t_k})\nonumber\\
	&\leq \max_{\mathbf{d}_{t_k,N+1},\mathbf{v}_{t_k,N+1}} J_{N+1}^1(x_{t_k},\tilde{\mathbf{u}}_{t_k,N+1},\mathbf{v}_{t_k,N+1},\mathbf{d}_{t_k,N+1},\beta_{t_k})\nonumber\\
	&\leq \max_{\mathbf{d}_{t_k,N},\mathbf{v}_{t_k,N}} J_{N}^1(x_{t_k},\check{\mathbf{u}}_{t_k,N},\mathbf{v}_{t_k,N},\mathbf{d}_{t_k,N},\beta_{t_k})+\bar{d}\nonumber\\
	&\leq V_N^1(x_{t_k},\beta_{t_k})+\bar{d},
	\end{align}
	where $\bar{d} = \max_{d\in\mathcal{D}}\alpha(\|d\|)$. Consequently, for all $x_{t_k}\in\mathcal{X}_f$, we have 
	\begin{align*}
	V_{N+1}^1(x_{t_k},\beta_{t_k})&\leq V_{1}^1(x_{t_k},\beta_{t_k}) + (N-1)\max_d\alpha(\|d\|)\\
	&\leq V_{1}(x_{t_k}) + \frac{1-\beta_{t_k}}{\beta_{t_k}}\ell(x_{t_k},\kappa_f(x_{t_k})) +\\
	&~~~~~~~(N-1)\bar{d}\\
	&\leq \ell_f(x_{t_k}) + \frac{1-\beta_{t_k}}{\beta_{t_k}}\ell(x_{t_k},\kappa_f(x_{t_k})) +N\bar{d}\\
	&\leq\sigma_3\|x_{t_k}\|+N\bar{d}.
	\end{align*}
	Hence, by induction, it is derived from the triggering condition (\ref{STAMPC:opt:self_trig_cond}) that $	V_N^{H_{t_k}^*}(x_{t_k},\beta_{t_k})\leq\sigma_3\|x_{t_k}\|+N\bar{d}, \forall x_{t_k}\in\mathcal{X}_f.$	For $x_{t_k}\in\mathcal{X}_N$ but $x_{t_k}\not\in\mathcal{X}_f$, the upper bound can be found by following Lemma 1 in \cite{limon_input_2006Auto}. Let $\mathcal{B}_r$ denote a set and  $\mathcal{B}_r = \{x\in\mathbb{R}^{n_x}| \|x\|\leq r\}$. Since $\mathcal{X},\mathcal{U},\mathcal{V}$ and $\mathcal{D}$ are compact sets, there definitely exists a finite constant $\bar{V}_N\in\mathbb{R}$ such that $	V_N^{H_{t_k}^*}(x_{t_k},\beta_{t_k})\leq \bar{V}_N$ for all $x_{t_k}\in\mathcal{X}_N$. For $x_{t_k}\in\mathcal{X}_N$ but $x_{t_k}\not\in\mathcal{B}_r$, we have $\|x_{t_k}\|\geq r$ and $	V_N^{H_{t_k}^*}(x_{t_k},\beta_{t_k})\leq \bar{V}_N$, which in turn leads to $	V_N^{H_{t_k}^*}(x_{t_k},\beta_{t_k})\leq \frac{\bar{V}_N}{r}\|x_{t_k}\|$. Consequently, we have  $V_{N}^{H_{t_k}^*}(x_{t_k},\beta_{t_k})\leq\sigma_{\Delta}\|x_{t_k}\|+N\bar{d}$	for all $x_{t_k}\in\mathcal{X}_n$, where $\sigma_{\Delta} = \max(\sigma_3,\frac{\bar{V}_N}{r}\|x_{t_k}\|)$.
	
	According to the triggering condition (\ref{STAMPC:opt:self_trig_cond}) and $\mathcal{V}_{l+t_k}\subseteq\hat{\mathcal{V}}_{l|t_k}$, for all $x_{t_k}\in\mathcal{X}_N$, we have 
	
	\begin{subequations}\label{STAMPC:V_k_V_k+1}
		\begin{align}
		&V_{N}^{H^*(x_{t_{k+1}})}(x_{k+1},\beta_{t_{k+1}}) - V_{N}^{H_{t_k}^*}(x_{t_k},\beta_{t_k})\nonumber\\
		\leq &V_{N}^{1}(x_{k+1},\beta_{t_{k+1}}) - V_{N}^{H_{t_k}^*}(x_{t_k},\beta_{t_k})\\
		\leq &V_{N}^{1}(x_{k+1},\beta_{t_{k+1}}) - \max_{v_{l|t_k}\in\hat{\mathcal{V}}_{l|t_k}, d_{l|t_k}\in\mathcal{D},l\in\mathbb{N}_{[0,H-1]} }\{\sum_{l=0}^{H_{t_k}^*-1}\nonumber\\
		&~~~~~\frac{1}{\beta_{t_k}}\ell(x_{l|t_k},u^*_{l|t_k})+V_{N-H_{t_k}^*}(x_{H|t_k})\}\\
		\leq & V_{N}^{1}(x_{k+1},\beta_{t_{k+1}}) - V_{N-H_{t_k}^*}(x_{t_{k+1}})\nonumber\\
		&~~~~~-\sum_{l=0}^{H_{t_k}^*-1}\frac{1}{\beta_{t_k}}\ell(x_{t_k+l},u^*_{l|t_k}).
		\end{align}	
	\end{subequations}
	
	It is similar to (\ref{STAMPC:V_N_Mono}) that $V_{N+1}(x_{t_{k+1}}) - V_{N}(x_{t_{k+1}})\leq \bar{d},$ which yields $V_{N}(x_{t_{k+1}}) - V_{N-H_{t_k}^*}(x_{t_{k+1}}) \leq H_{t_k}^*\bar{d}$ for $x_{t_{k+1}}\in\mathcal{X}_{N-H_{t_k}^*}$. 
	 Hence, we have $V_{N}^{H^*_{t_{k+1}}}(x_{k+1},\beta_{t_{k+1}}) - V_{N}^{H_{t_k}^*}(x_{t_k},\beta_{t_k})\leq  H_{t_k}^*\bar{d} - \sum_{l=0}^{H_{t_k}^*-1}\frac{1}{\beta_{t_k}}\ell(x_{t_k+l},u^*_{l|t_k}) \leq H_{t_k}^*\bar{d} - \sum_{l=0}^{H_{t_k}^*-1}\frac{1}{\beta_{t_k}}\ell(x_{t_k+l},u^*_{l|t_k}).$ In summary, it is inferred that, by applying the proposed self-triggered adaptive MPC scheme presented in Algorithm \ref{STAMPC:alg:stampc_state_para_est} to the system in (\ref{STAMPC:sys_nonlinear_affine}), the corresponding optimal value function $V_{N}^{H_{t_k}^*}(x_{t_k},\beta_{t_k})$ is an ISpS Lyapunov function. According to \cite[Theorem 1]{limon_input_2006Auto}, it is proved that the closed-loop system in (\ref{STAMPC:sys_nonlinear_affine}) is ISpS under the proposed self-triggered adaptive MPC algorithms for all $x_0\in\mathcal{X}_N$ at triggering time instants.
\end{IEEEproof}
\begin{remark}
Compared with the self-triggered robust method, the time-varying $\beta_k$ is considered and updated with respect to $\hat{\mathcal{V}}_{t_k}$ in the proposed method. This strategy improves the effective use of the EFSS in MPC framework to reduce the conservatism. The design of $\beta_{t_k}$ in (\ref{STAMPC:beta_k}) is to further enlarge the triggering interval when a less conservative EFSS is obtained. For a better control performance, we can choose a smaller $\beta_k$ if the size of $\hat{\mathcal{V}}_{t_k}$ becomes smaller. Therefore, a suitable evolution of $\beta_k$ helps achieve a trade-off between the communicational load and performance. 
\end{remark}
% % % % % %% % % % % % % % % % % % % % % % % % % % % % % % % % % % % % % % % % % % % % % % % % % % % % % % % % % % % % % % % % % % % % % % % % % % % % % %
\section{Illustrative example}\label{STAMPC:sec:simulation}
In this section, a numerical example is presented to validate our theoretical results. Consider the following discrete-time cart and spring-damper system from \cite{acipl_liu_robust_2019Auto,raimondo_min_2009EJC}
\begin{align*}
x_{t+1}(1) = &x_t(1)+Tx_t(2),\\
x_{t+1}(2) = &-\frac{lT}{m}e^{-x_t(1)}+\frac{m-hT}{m}x_t(2)+\frac{T}{m}u_t\\
&-\frac{T}{m}v_tx_t(2)+\frac{T}{m}d_t.
\end{align*}
where $x_t=[x_t(1)~x_t(2)]^\T$ and $u_t$ are the system state and input satisfying the constraints $|u_t|\leq 4.5~N$ and $|x_t(1)|\leq 2~m$. The additive and parametric uncertainties are limited by $|d_t|\leq 0.1, |v_t|\leq 0.15$ and $\|v_{t+1}-v_t\|\leq 0.008$. More specifically, we consider the following sequence of uncertainties in this example: $v_t = 0.1\sin(\frac{4t}{\pi})$ and $d_t = 0.15\sin(\frac{t}{\pi})$.
The system parameters are given as follows: $m=1~kg;~l=0.33~N/m;~h=1.1~Ns/m;~T=0.4~s$. For the proposed self-triggered adaptive MPC scheme, we set the prediction horizon $N = 6$. The stage cost function is chosen as $\ell(x,u) = x^\T Q x+u^\T R u$ with $Q = \text{diag}(0.64,0.64)$ and $R = 1$. According to \cite{raimondo_min_2009EJC,acipl_liu_robust_2019Auto}, we set $\mathcal{X}_f = \{x:x^\T P x\leq 3.8 \}$ with $P = \begin{bmatrix}4.5678&3.2018\\3.2018&4.3500\end{bmatrix}$. $\ell_f(x) = x^\T P x$; $\kappa_f(x) = [-0.7797~-1.1029]x$.  The feedback policy in (\ref{STAMPC:opt_closed}) is set as $u(x) = a\kappa_f(x)+bx^\T x+c$ where $a,b,c\in\mathbb{R}$ are the decision variables for the optimization problem in (\ref{STAMPC:opt_closed}). To demonstrate the effectiveness of the proposed methods, the self-triggered robust MPC \cite{acipl_liu_robust_2019Auto} (ST-RMPC) is also implemented with the same parameters for the purpose of comparison. For the self-triggering scheduler, we set $H_{\max} = 5$ and $\beta_0 = 1.1$. 

Starting from the initial condition $x_0 = [1~1]^\T$, the trajectories of system states are shown in Fig. \ref{STAMPC:fig:trajx1x2}, and the time evolution of control inputs is plotted in Fig. \ref{STAMPC:fig:traju}. The triggering time instants are reported in Fig. \ref{STAMPC:fig:triginst}. To further illustrate the performance of the proposed methods, we introduce the regulation performance index $J_p = \sum_{t=0}x_t^\T Qx_t+u_t^\T Ru_t$. Table \ref{STAMPC:tab:per_com} shows the comparison of performance index and average sampling time, where we assume that the average sampling time for the periodic MPC is 1.  It can be seen from Figs. \ref{STAMPC:fig:trajx1x2}-\ref{STAMPC:fig:traju} and Table \ref{STAMPC:tab:per_com} that both state and input constraints are satisfied for all $t\geq 0$, but the proposed self-triggered adaptive MPC method can further reduce the average sampling frequency,  while the comparable closed-loop performance is guaranteed. Fig. \ref{STAMPC:fig:trajset} demonstrates the evolution of EFSS for unknown parameters, where the red mark indicates real value of $v_t$ and the length of the blue vertical line indicates the size of EFSS. It is worthwhile to observe that the falsified parameter can be removed by using the proposed set-membership parameter estimator under the self-triggering mechanism.
%\begin{figure}
%	\centering
%	\includegraphics[width=0.9\linewidth]{traj_x1}
%	\caption{The time evolution of system state $x_t(1)$.}
%	\label{STAMPC:fig:trajx1x2}
%\end{figure}
%
%\begin{figure}
%	\centering
%	\includegraphics[width=0.9\linewidth]{traj_x2}
%	\caption{The time evolution of system state $x_t(2)$.}
%	\label{STAMPC:fig:trajx2}
%\end{figure}

\begin{figure}
	\centering
	\includegraphics[width=0.82\linewidth]{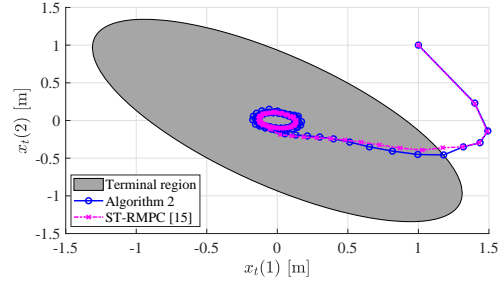}
	\caption{The time evolution of system states.}
	\label{STAMPC:fig:trajx1x2}
\end{figure}
\begin{figure}
	\centering
	\includegraphics[width=0.82\linewidth]{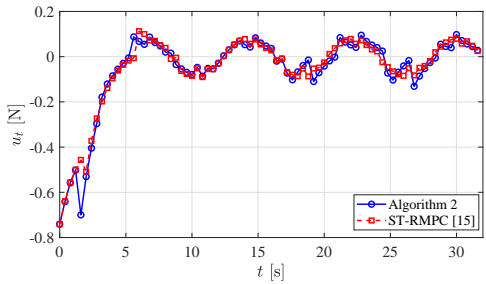}
	\caption{The time evolution of control input $u_t$.}
	\label{STAMPC:fig:traju}
\end{figure}

\begin{figure}
	\centering
	\includegraphics[width=0.82\linewidth]{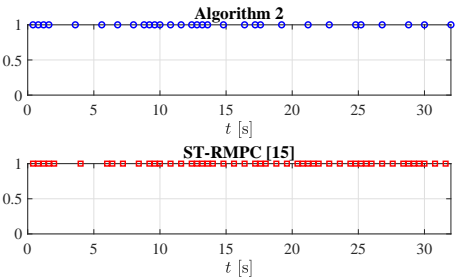}
	\caption{Triggering time instants for the proposed method and the ST-RMPC.}
	\label{STAMPC:fig:triginst}
\end{figure}

%the lower computation frequency can be achieved by improving the accuracy of the set-membership estimation. 
\begin{table}[!t]
	\renewcommand{\arraystretch}{1.3}
	\centering
	%\caption{Comparison of regulation performance and average sampling time}
	\caption{Closed-loop performance comparison.}
	\begin{tabular}{ccc}
		\hline
		&\bfseries $J_p$ & \bfseries Average sampling time\\
		\hline
		% Algorithm \ref{STAMPC:alg:stampc_para_est}&	7.3422 & 2.7949\\
		Algorithm \ref{STAMPC:alg:stampc_state_para_est}&	13.4122 & 2.6333\\
		ST-RMPC\cite{acipl_liu_robust_2019Auto}&	13.1229 & 1.6250\\
		\hline
	\end{tabular}
	\label{STAMPC:tab:per_com}
\end{table}

\begin{figure}
	\centering
	\includegraphics[width=0.82\linewidth]{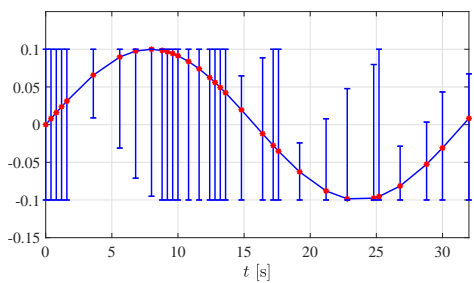}
	\caption{Feasible solution set for the time-varying unknown parameters estimated at triggering time instants.}
	\label{STAMPC:fig:trajset}
\end{figure}

% % % % % % % % % % % % % % % % % % % % % % % % % % % % % % % % % % % % % % % % % % % % % % % % % % % % % % % % % % % % % % % % % % % % % % % % % % % % % %
\section{Conclusion}\label{STAMPC:sec:conclusion}
In this work, we developed a self-triggered adaptive  MPC approach for constrained discrete-time nonlinear systems subject to parametric uncertainties and additive disturbances. A set-membership parameter estimator was designed to estimate the FSS of unknown parameters by using the indirect polytopic set computation. The zonotopic description of the FSS was considered to reduce overestimation. By integrating the proposed estimator with the min-max MPC technique, the robust constraint satisfaction was guaranteed with reduced conservatism. The new estimation of uncertainties was considered in both the triggering scheduler and the formulation of MPC optimization problem, therefore giving rise to the decreased sampling frequency compared with the robust self-triggered method. It has been proved that the proposed self-triggered adaptive MPC algorithm is recursively feasible and the closed-loop system is ISpS. We illustrated the efficacy and advantages of the proposed method over \cite{acipl_liu_robust_2019Auto} by using a numerical example.    
% % % % % % % % % % % % % % % % % % % % % % % % % % % % % % % % % % % % % % % % % % % % % % % % % % % % % % % % % % % % % % % % % % % % % % % % % % % % % %
\bibliographystyle{ieeetran}
\bibliography{ref}
\end{document}